\documentclass{amsart}
\usepackage{amssymb,amstext,amsthm,titletoc,fancyhdr,color,ulem}

\newtheorem{theorem}{Theorem}[section]
\newtheorem{lemma}[theorem]{Lemma}

\theoremstyle{definition}

\theoremstyle{remark}
\newtheorem{remark}[theorem]{Remark}

\numberwithin{equation}{section}

\begin{document}

\title[On the growth of $\Delta^n f(z)/f(z)$ and Wiman-Valiron theory] {On the growth of 
 logarithmic difference of meromorphic functions and a Wiman-Valiron estimate}

\author{Yik-Man Chiang}
\address{Department of Mathematics, Hong Kong University of Science and
Technology, Clear Water Bay, Kowloon, Hong Kong,\ P. R. China.}
\email{machiang@ust.hk}
\thanks{The first and second authors were partially supported by the Research Grants Council of the Hong Kong Special Administrative Region, China (GRF600609 and GRF601111) and the HKUST PDF Matching Fund. The second author was also partially supported by the National Natural Science Foundation of China  (Grant No. 11271352).}

\author{Shao-Ji Feng}
\address{Academy of Mathematics and Systems Science, Chinese Academy
of Sciences, Beijing, 100080, P. R. China.}
\email{fsj@amss.ac.cn}


\date{18th December 2015
}

\dedicatory{Dedicated to the memory of Milne Anderson.}


\begin{abstract}
The paper gives a precise asymptotic relation between higher order logarithmic difference and logarithmic derivatives for meromorphic functions with order strictly less then one. This allows us to formulate a useful Wiman-Valiron type estimate for logarithmic difference of meromorphic functions of small order. We then apply this estimate to prove a classical analogue of Valiron about entire solutions to linear differential equations with polynomials coefficients for linear difference equations.
\end{abstract}

\maketitle



\section{Introduction}
\label{intro}
It was shown by Ablowitz, Herbst and Halburd \cite{AHH} and Halburd and Korhonen \cite{HK-2} (see also \cite{HK-5}) that the integrability of the discrete Painlev\'e equations in the complex plane can be characterised by the finite Nevanlinna order of growth of meromorphic solutions. A crucial role is played in their Nevanlinna theory approach is some difference versions of the logarithmic derivative estimates which were established in \cite{HK-1} and independently by us in \cite{CF-1}. 
Such estimates were applied in a number of applications recently (see e.g. \cite{HKLRT}).
 We show in a subsequent paper \cite{CF-2} that given $\varepsilon>0$, there exists a set $E\subset (1,\infty)$ of $|z|=r$ of finite logarithmic measure, so that
\begin{equation}
\label{E:Fund-thm-est-2}
    \frac{f(z+\eta)}{f(z)}=e^{\eta\frac{f^\prime(z)}{f(z)}+
    O(r^{\beta+\varepsilon})},
\end{equation}
holds for $r\not\in E\cup[0,1]$, where 
    $\beta=\max\{\sigma-2,\,2\lambda-2\}$ if $\lambda<1$ and
    $\beta=\max\{\sigma-2,\,\lambda-1\}$ if $\lambda\ge1$ and
    $\lambda=\max\{\lambda^\prime,\lambda^{\prime\prime}\}$.
    We deduce from the equation \eqref{E:Fund-thm-est-2} that
when the order $\sigma(f)<1$, then for each given $\varepsilon>0$, there is
a set $E\subset (1,\, \infty)$ of finite logarithmic measure such that 
\begin{equation}
\label{E:delta-deriv-1} \frac{\Delta^nf(z)}{f(z)}=\eta^n
\frac{f^{(n)}(z)}{f(z)}+O(r^{(n+1)(\sigma-1)+\varepsilon})
\end{equation}
holds for $|z|=r\notin E$, where $\Delta f(z)=f(z+\eta)-f(z)$, $\Delta^nf(z)=\Delta\big(\Delta^{n-1} f\big)$. The above results show the different behaviours for meromorphic functions of order less than and greater than unity. Bergweiler and Langley \cite[Lem. 4.2]{BL2007} also obtained an asymptotic behaviour $\Delta^n f(z)\sim f^{(n)}(z)$ as $z\to\infty$ outside an $\varepsilon-$set. However no precise error bounds were given.
\medskip

If one adopts the symbolic operator $D=\frac{d}{dz}$, then it was discussed in \cite[\S 8]{CF-2} that we could write down the formal expression
\smallskip

\begin{equation*}
\label{E:taylor-expan-1}
    \Delta^nf =\big(\eta^n D^n+ \frac{n}{2!}\eta^{n+1} D^{n+1}+\cdots\big)f.
\end{equation*}
\smallskip
This note has two purposes. We first establish rigorously that the above  formal expression is indeed valid for meromorphic functions of order strictly less than one. This allows us to obtain a more precise difference Wiman-Valiron estimate than that established in \cite{CF-2}. 
\bigskip

Before we can state our first main result, let us consider the well-known formal expansion

\begin{equation}\label{E:Formal-exp} 
	\Delta^nf(z)=n!\, \sum_{k=n}^\infty \frac{\mathfrak{S}_k^{(n)}}{k!}\,f^{(k)}(z),
\end{equation}
where the $\mathfrak{S}_k^{(n)}$ are \textit{Stirling numbers of the second kind} \cite[\S 24.1]{AS64}. We recall that Stirling number of the second kind $\frak{S}_n^{(m)}$ counts the number of different ways to partition a set of $n$ objects into $m$ non-empty subsets. In particular, it has the following generating function \cite[\S21.1.4]{AS64}:
	\[
		x^n=\sum_{m=0}^n \frak{S}_n^{(m)}\, x(x-1)\cdots (x-m+1).
	\]
\bigskip

We now give a rigorous justification of the formal expansion (\ref{E:Formal-exp}).
\bigskip
	
\begin{theorem}  \label{T:higer-difference} Let $f$ be a meromorphic function with order $\sigma=\sigma(f)<1$. Then for any positive integers $n,\, N$ such that $N\ge n$, and for each $\varepsilon>0$,  there is a set $E\subset [1,\, +\infty)$ of finite logarithmic measure so that
\begin{equation}\label{E:higher-difference}
	\frac{\Delta^n f(z)}{f(z)}=n!\, \Bigg(\sum_{k=n}^N\frac{\eta^k\mathfrak{S}_k^{(n)}}{k!}\frac{f^{(k)}(z)}{f(z)}\Bigg)+O\big(\eta^{n+N+1}r^{(n+N+1)(\sigma-1)+\varepsilon}\big)
\end{equation}
for $|z|=r\notin E\cup [0,1]$
\end{theorem}
\bigskip

A classical result states that entire solutions to linear differential equations with polynomial coefficients have {completely regular growth} of rational order (see Valiron \cite[pp. 106-108]{Val}). The second purpose of this paper is to prove an analogue to this result that entire solutions to linear difference equations with polynomial coefficients of order strictly less than one must have \textit{completely regular growth} of rational order. Previously, Ishizaki and Yanagihara \cite{IY2004} established the same result for entire solutions of order strictly less than $1/2$ via a different method. In \cite[Theorem 7.3]{CF-2} the authors  extended the theorem of Ishizaki and Yanagihara to include transcendental entire solutions of order strictly less than unity must have rational order of growth, but they failed to establish that the growth of these solutions must be completely regular.
\bigskip

This paper is organised as follows. Some preliminary results that are needed for the proof of Theorem \ref{T:higer-difference} are given in \S\ref{S:lemmas}, followed by the proof of Theorem \ref{T:higer-difference} in \S\ref{S:proof}. We then formulate sharp difference-type Wiman-Valiron estimates in \S\ref{S:Wiman-Valiron} and apply them to prove a difference version of a result originally for differential equations  stated  in Valiron \cite{Val} in the next section \S\ref{S:application}, which is the second main purpose of this paper.

\section{Lemmas}
\label{S:lemmas}

\begin{theorem} \label{T:first-difference} Let $f$ be a meromorphic function with order $\sigma=\sigma(f)<1$. Then for each positive integer $k$, and for each $\varepsilon>0$, 
there exists an exceptional set $E^{(\eta)}$ in $\mathbb{C}$ consisting of a union of disks centred at the zeros and poles of $f(z)$  such that when $z$ lies outside of the $E^{(\eta)}$,
	\begin{equation}\label{E:first-difference}
		\begin{split}
	\frac{\Delta f}{f}: &=\frac{f(z+\eta)-f(z)}{f(z)}\\
						&=\eta \frac{f^\prime(z)}{f(z)} +\frac{\eta^2}{2!}\frac{f^{\prime\prime}(z)}{f(z)}+\cdots+ \frac{\eta^k}{k!}\frac{f^{(k)}(z)}{f(z)}+O\big(\eta^{k+1}r^{(k+1)(\sigma-1)+\varepsilon}\big).
		\end{split}
	\end{equation}
Moreover, the set $\pi E^{(\eta)}\cap [1,\, +\infty)$, where $\pi E^{(\eta)}$ is obtained from rotating the exceptional disks of $E^{(\eta)}$ so that their centres all lie on the positive real axis, has finite logarithmic measure.
\end{theorem}
\bigskip

We recall that a subset $E$ of $\mathbb{R}$ has \textit{finite logarithmic measure} if
\[
	\textrm{lm}(E)=\int_{E\cap (1,\,\infty)} \frac{dr}{r}
\]
is finite. The set $E$ is said to have an infinite logarithmic measure if lm$(E)=+\infty$.
\medskip

We need a \textit{uniformity} estimate for the logarithmic derivatives and logarithmic difference estimates that hold outside an exceptional set of $|z|=r$ of finite logarithmic measure \textit{simultaneously}. So let us first  review Gundersen's important logarithmic derivative estimate \cite[Cor. 2]{Gund88} and our earlier difference analogous estimate \cite[Thm. 8.2]{CF-1}.
\bigskip

\begin{lemma}[\cite{Gund88} Lem. 8]
\label{L:Gundersen}
	Let $f$ be a transcendental meromorphic function, and let $k$ be a positive integer. Let $\{c_n\}$ be the sequence of zeros and poles of $f$,  where $\{c_n\}$ is listed according to multiplicity and ordered by increasing modulus. Let $\gamma>1$ be a given real constant. Then there exist constants $R_0>0$ and $C=C(\gamma,\, k)>0$ such that for $|z|\ge R_0$ and $f(z)\not=0,\, \infty$, we have
\begin{equation}\label{E:Gundersen}
	\Big|\frac{f^{(k)}(z)}{f(z)}\Big|\le C\Big(\frac{T(\gamma r,\, f)}{r}+
	\sum_{|c_k|<\gamma r} \frac{1}{|z-c_k|}\Big)^k,
\end{equation}
where $r=|z|$.
\end{lemma}
\bigskip

\begin{lemma}[\cite{CF-1} Eqn. (8.5)]\label{L:log-difference}
	Let $f$ be a transcendental meromorphic function, and let $\{c_n\}$ be the sequence of zeros and poles of $f$,  where $\{c_n\}$ is listed according to multiplicity and ordered by increasing modulus. Let
$\gamma>1$ and a complex $\eta$ be given, then there are constants $R_1>\frac{2}{\gamma}\,|\eta|>
0$, $\beta=\beta(\gamma)$ such that for $|z|\ge R_1$,  we have
\begin{equation}\label{E:direct-estimate}
	\bigg|\log\bigg|\frac{f(z+\eta)}{f(z)}\bigg|\bigg| \leq
	|\eta|\,\beta\,  \frac{T\big(\gamma r,\, f\big)}{r}
    +|\eta|\cdot\sum_{|d_k|<\gamma r}
    \frac{1}{|z-d_k|}
\end{equation}
holds, where $|z|=r$ and  $(d_k)_{k\in N}=(c_k)_{k\in N}\cup(c_k-\eta)_{k\in N}$.
\end{lemma}
\bigskip

\begin{lemma}
\label{L:thm-ptwse-2} Let $f(z)$ be a meromorphic function of finite order
$\sigma$, $\eta$ a non-zero complex number, and  $\varepsilon>0$
be a given real constant. Then there exists an exceptional set $E^{(3\eta)}$ in $\mathbb{C}$ consisting of a union of disks centred at the zeros and poles of $f(z)$, as described in the Theorem \ref{T:first-difference}, such that when $z$ lies entirely outside the $E^{(3\eta)}$, then 
 we have
\begin{enumerate}
	\item[]
		\begin{equation}
\label{E:thm-Gund-b}
 \bigg|\frac{f^{(k)}(z)}{f(z)}\bigg|\le
    |z|^{k(\sigma-1)+\varepsilon},
\end{equation}
and \item[]
\begin{equation}
\label{E:thm-ptwse-2-b}
    \exp\big(-|z|^{\sigma-1+\varepsilon}\big)\leq
    \bigg|\frac{f(z+t\eta)}{f(z)}\bigg|\leq
    \exp\big(|z|^{\sigma-1+\varepsilon}\big)
\end{equation}
\end{enumerate}
hold simultaneously and uniformly in $t\in [0,\,1]$. Moreover, the set $\pi E^{(3\eta)}\cap (0,\, +\infty)$, where $\pi E^{(3\eta)}$ is obtained from rotating the exceptional disks of $E^{(3\eta)}$ so that their centres all lie on the positive real axis, has 
finite logarithmic measure. 
\end{lemma}
\bigskip

The above inequality (\ref{E:thm-Gund-b}) can be derived from (\ref{E:Gundersen}) in \cite[Lem. 8]{Gund88}, while the inequality (\ref{E:thm-ptwse-2-b}) from (\ref{E:direct-estimate})  in \cite[Eqn. (8.5)]{CF-1}. Both inequalities hold generally under the finite order assumption on $f$, and both proofs require the  Poisson-Jensen formula and the classical Cartan lemma, except that we need to take into the account of the \textit{uniformity assumption} on $t$. However, in order to fulfill our later application to establishing the Theorem \ref{T:first-difference}, one needs more detailed information concerning the exceptional set that arises from removing the zeros and poles from applying the Cartan lemma (below) than quoting Gundersen \cite[\S 7]{Gund88} directly. So we judge it is appropriate to offer a full proof of the Lemma \ref{L:thm-ptwse-2} based on that in \cite{Gund88} but tailored to our need here. In fact, our construction of the exceptional disks which have larger radii than those constructed in \cite{Gund88} by $3\eta$. This is to guarantee the fact that we need below, namely that whenever $z$ lies outside $E^{(3\eta)}$, then the line segment $[z,\, z+\eta]$ lies outside $E^{(\eta)}$.

\bigskip

Let us first recall Cartan's theorem \cite{Cartan1928} which we adopt from \cite{Gund88}.

\begin{lemma}[\cite{Cartan1928}]\label{L:Cartan} Let $a_1,\, \cdots, a_m$ be any finite collection of complex numbers, and let $d>0$ be any given positive number. Then there exists a finite collection of closed disks $D(a_k,\, r_k)\ (1\le k\le m)$ with corresponding radii that satisfy\\ $\sum_{k=1}^mr _k=2d$, and a permutation of the points $a_1,\, \cdots, a_m$, labeled by $b_1,\, \cdots, b_m$, say, such that for $z\not\in \bigcup_{k=1}^mD(a_k,\, r_k)$,
	\begin{equation}
		\label{E:Cantan}
			|z-b_k|>\frac{k}{m}\, d, \quad 1\le k\le m,
	\end{equation}
where the permutation may depend on $z$.
\end{lemma}

\subsubsection*{Proof of Lemma \ref{L:thm-ptwse-2}}
\begin{proof} Let us first choose $\gamma >1$ and let us define the annulus

	\begin{equation}
		\label{E:annulus}
			\Gamma_\nu:=\{z:\, \gamma^\nu\le |z|\le \gamma^{\nu+1}\},\quad \nu\in\mathbb{N}.
	\end{equation}

We  let $R=\gamma^{\nu+2}$, $d_k(t)=(c_k)\cup (c_k-t\eta),\ (0\le t\le 1)$ so that $(d_k(0))=~(c_k)$ and $(d_k(1))=(d_k)$. Here the $d_k,\, c_k$ are defined in the Lemma \ref{L:log-difference}. Let $m=n(R)=n(\gamma^{\nu+2})$  
so that $|d_m(1)|\le R$ and $|d_{m+1}(1)|>R$.   Let us now choose an integer $\nu_0$ such that
	\begin{equation}
		\label{E:nu-1st}
			\gamma^{\nu_0+2}\ge |d_1(1)|
	\end{equation}
and
	\begin{equation}
		\label{E:nu-2nd}
			1\le (\gamma-1)\,\log n(\gamma^{\nu_0+2}).
	\end{equation}
We choose
	\begin{equation}
		\label{E:diam_3eta}
			d_{3\eta}:=\frac{\gamma^\nu}{[\log (\gamma^\nu)]^\gamma}+3|\eta|.
	\end{equation}
Lemma \ref{L:Cartan} asserts that there is a finite collection of closed disks $D(a_k,\, r_k)\ (1\le k\le m)$ whose radii have a sum of $2d_{3\eta}$ and that when $z\not\in \bigcup_{k=1}^mD(a_k,\, r_k)$, then there is a permutation of the points $a_1,\, \cdots, a_m$, re-labelled by $b_1,\, \cdots, b_m$, say,  that, according to \eqref{E:Cantan},  satisfies
\medskip

	\begin{equation}
		\label{E:Cantan_applied}
			|z-b_k|>\frac{k}{m}\, \Big(\frac{\gamma^\nu}{[\log (\gamma^\nu)]^\gamma}+3|\eta|\Big), \quad 1\le k\le m.
	\end{equation}

Let us confine our $z$ within the annulus (\ref{E:annulus}). Because of the choice of $\nu_0$ as defined in (\ref{E:nu-1st}) and \eqref{E:nu-2nd}, we deduce 
	\begin{align}
		\label{E:upper-bound}
				\sum_{|d_k|<\gamma r}   \frac{1}{|z-d_k(1)|} &\le 
				\sum_{k=1}^m \frac{1}{|z-b_k|}\le \frac{m [\log (\gamma^\nu)]^\gamma}{\gamma^\nu+3|\eta|}\sum_{k=1}^m\frac1k\notag \\
				&\le  \frac{n(R) [\log (\gamma^\nu)]^\gamma}{\gamma^\nu+3|\eta|} (1+\log m)\notag \\
				&\le \frac{n(\gamma^2 r)\log^\gamma r}{\gamma^\nu+3|\eta|}\big(1+\log n(\gamma^2r)\big)\notag \\
				&< \frac{n(\gamma^2 r)\log^\gamma r}{\gamma^\nu}\big(1+\log n(\gamma^2r)\big)\notag \\
				&\le \gamma\frac{n(\gamma^2 r)}{r}\log^\gamma r\big(1+\log n(\gamma^2r)\big)\notag \\
				&\le \gamma^2\frac{n(\gamma^2 r)}{r}\log^\gamma r\log n(\gamma^2r).
	\end{align}
\bigskip
So given an $\varepsilon>0$, we deduce from the inequalities (\ref{E:Gundersen}) and (\ref{E:direct-estimate}) and (\ref{E:upper-bound}) that there is an exceptional set $ E^{(3\eta)}=\bigcup_{\nu=1}^\infty E^{(3\eta)}_\nu$ where each $E^{(3\eta)}_\nu$ consists of the union of closed disks $D(a_k,\, r_k)\ (1\le k\le m)$ whose centres lie in $\Gamma_\nu$  defined above
 such that for all $z\not\in E^{(3\eta)}$, the (\ref{E:thm-Gund-b}) and (\ref{E:thm-ptwse-2-b}) hold simultaneously and uniformly for all $t\in [0,\, 1]$. So we can choose an $R$ so large such that, when $|z|>R$ and $z\not\in E^{(3\eta)}$,
the inequalities (\ref{E:thm-Gund-b}), (\ref{E:thm-ptwse-2-b}) and (\ref{E:upper-bound})
 hold simultaneously.
\bigskip

It remains to compute the size of the exceptional sets. According to Lemma \ref{L:Cartan} and \eqref{E:diam_3eta} that the sum of the diameters of the disks $\bigcup_{k=1}^m D(a_k,\, r_k)$ is
	\begin{equation}
		\label{E:4diam_3eta}
			4\, d_{3\eta}=4\Big(\frac{\gamma^\nu}{[\log (\gamma^\nu)]^\gamma}+3|\eta|\Big).
	\end{equation}

We argue in the spirit of \cite{Gund88} that we revolve each of the disks $D(a_k,\, r_k)$ about the origin to form annuli centered about the origin. We then consider the logarithmic measure of the union of the intersection of these annuli with the positive real axis in $\Gamma_\nu$ which we denote by $\pi E^{(3\eta)}$. Here we need to clarify an exceptional situation, namely when one of those disks already contains the origin, then we simply count the line segment $[0,\, |a_k|]$ part of the exceptional set in the computation of the logarithmic measure. That is, we have
	\begin{equation*}
		E_\nu := [\gamma^\nu,\, \gamma^{\nu+1}] \cap \pi E^{(3\eta)}
	\end{equation*}
and
	\begin{equation}
		\label{E:exceptionl}
			 E:=\bigcup_{\nu= \nu_0}^\infty E_\nu.
	\end{equation}
\medskip

It follows from (\ref{E:4diam_3eta}) that the linear measure of $E_\nu$ does not exceed $4\, d_{3\eta}$. So 
	\begin{equation*}
		\begin{split}
		\int_{E}\frac{dx}{x} &=\sum_{\nu=\nu_0}^\infty \int_{E_\nu}\frac{dx}{x}
		\le \sum_{\nu=\nu_0}^\infty
			\Big\{\log \Big(\gamma^\nu+\frac{4\gamma^\nu}{[\log (\gamma^\nu)]^\gamma}+12|\eta|\Big)
			-\log \gamma^\nu\Big\}\\
			&= \sum_{\nu=\nu_0}^\infty \log \Big( 1+\frac{4}{[\log (\gamma^\nu)]^\gamma}
			+\frac{12|\eta|}{\gamma^\nu}\Big)<\infty
		\end{split}
	\end{equation*}
since $\gamma>1$. This proves that the exceptional set (\ref{E:exceptionl}) has finite logarithmic measure.
\end{proof}
\bigskip

\begin{remark} We would like to emphasis that we take the exceptional set to be $ E^{(\eta)}$  in our application below. Since the  $ E^{(\eta)}$ is a subset of $ E^{(3\eta)}$, so the set $E^{(\eta)}$ also has finite logarithmic measure.
\end{remark}
\bigskip

We will also require a complex form of Lagrange's version of Taylor's theorem. Although we cannot find an exact reference for the result, one can easily modify the argument in \cite[p. 242]{Apo} to obtain:
\medskip

\begin{lemma}\label{L:Taylor} Let $f$ be an analytic function in a domain $D$. Let $c\in D$, then
\begin{equation}\label{E:Taylor}
	f(z)=f(c)+f^\prime(c)(z-c)+\frac{f^{\prime\prime}(c)}{2!}(z-c)^2+\cdots+\frac{f^{(n)}(c)}{n!}(z-c)^n+R_n(z)
\end{equation}
where
\begin{equation}\label{E:reminder}
	R_n(z)=\frac{1}{n!}\int_c^z(z-t)^nf^{(n+1)}(t)\, dt
\end{equation}
for each $z\in D$.
\end{lemma}
\bigskip

\section{Proof of Theorem \ref{T:first-difference}}

We apply the Lemma  \ref{L:thm-ptwse-2} with $z\not\in E^{(3\eta)}$. Then the inequalities (\ref{E:thm-Gund-b}) and (\ref{E:thm-ptwse-2-b}) both hold. In fact, it is not difficult to see that the whole line segment $[z,\, z+\eta]$ lies outside $E^{(\eta)}$ which also has finite logarithmic measure. So let us choose a path in the complex plane that connects $z$ and $z+\eta$ and that does not intersect with the $z$ from the exceptional set $E^{(\eta)}$. We replace $z$ by $z+\eta$ and $c$ by $z$ in (\ref{E:Taylor}), and divide through both sides of the Taylor expansion (\ref{E:Taylor})  by $f(z)$ to yield
 \begin{equation*}\label{E:Taylor-2}
	\frac{f(z+\eta)-f(z)}{f(z)}=\eta\frac{f^\prime (z)}{f(z)}+\frac{\eta^2}{2!}\frac{f^{\prime\prime}(z)}{f(z)}+\cdots+\frac{\eta^n}{n!}\frac{f^{(n)}(z)}{f(z)}+\frac{R_n(z+\eta)}{f(z)}
\end{equation*}
where $R_n$ is the remainder given in \eqref{E:reminder}. Hence
\begin{align}\label{E:reminder-2}
	\frac{R_n(z+\eta)}{f(z)}&=\frac{1}{n!}\int_z^{z+\eta}(z+\eta-t)^n\frac{f^{(n+1)}(t)}{f(z)}\, dt\notag\\
		&=\frac{\eta^{n+1}}{n!}\int_0^1(1-T)^n\frac{f^{(n+1)}(z+T\eta)}{f(z)}\, dT.
\end{align}
\bigskip

Lemma \ref{L:thm-ptwse-2} asserts that for each $\varepsilon>0$, one can have an exceptional set $\pi E^{(\eta)}$ of real numbers such that (\ref{E:thm-Gund-b}) and (\ref{E:thm-ptwse-2-b}) hold simultaneous and uniformly for $t\in[0,\,1]$ and for $r=|z|$ outside of $\pi E^{(\eta)}$. Thus (\ref{E:reminder-2}) becomes
\begin{align}\label{E:Taylor-3}
	\Big|\frac{R_n(z+\eta)}{f(z)}\big| &=\Big|\frac{\eta^{n+1}}{n!}\int_0^1(1-T)^n\frac{f^{(n+1)}(z+T\eta)}{f(z+T\eta)}\,\frac{f(z+T\eta)}{f(z)}\, dT\Big|\notag\\
	& \le  \frac{|\eta|^{n+1}}{n!}\,2^n \cdot O\big(e^{|z+\eta T|^{\sigma-1+\varepsilon}}\big)\cdot O\big(|z+T\eta|^{(n+1)(\sigma-1)+\varepsilon}\big)\notag\\
	& \le O\big(|\eta|^{n+1} |z|^{(n+1)(\sigma-1)+\varepsilon}\big)
\end{align}
\bigskip

\noindent as $|z|\to\infty$ and outside of $\pi E^{(\eta)}$ which is a subset of 
$E$ from (\ref{E:exceptionl}). This proves (\ref{E:first-difference}).\qed
\medskip

\section{Proof of Theorem \ref{T:higer-difference}}\label{S:proof}

\begin{proof} 
We shall also make use of the recurrence formula \cite[p. 825]{AS64}
	\begin{equation}
		\label{E:recurrence}
	\binom{m}{r}\mathfrak{S}_n^{(m)}=\sum_{k=m-r}^{n-r}\binom{n}{k}
\mathfrak{S}_{n-k}^{(r)} \mathfrak{S}_{k}^{(m-r)} ,\quad n\ge m\ge r.
	\end{equation}
\medskip

We first note that if $n=1$, then because $\mathfrak{S}_k^{(1)}=1$ for all $k\in \mathbb{N}$ \cite[p. 825]{AS64} so the result (\ref{E:higher-difference}) follows from (\ref{E:first-difference}) of Theorem \ref{T:first-difference}. We apply induction on $n$. Let us suppose that the equation (\ref{E:higher-difference}) holds and let $F(z)=\Delta^n f(z)$. We now apply Theorem \ref{T:first-difference} to $F$ and the inductive hypothesis to obtain

\begin{align}
	\frac{\Delta^{n+1}f(z)}{f(z)} &=\frac{\Delta F(z)}{F(z)}
	\frac{F(z)}{f(z)}\notag\\
	&= \frac{F(z)}{f(z)} \Bigg(n!\Big(\sum_{k=n}^N\frac{\eta^k}{k!}\frac{F^{(k)}(z)}{F(z)}\Big)+O\big(\eta^{n+N+1}r^{(n+N+1)(\sigma-1)+\varepsilon}\big)\Bigg)\notag\\
	&= n!\sum_{k=1}^N\frac{\eta^k}{k!}\frac{\Delta^n (f(z))}{f(z)}\,\frac{[\Delta^n (f(z))]^{(k)}}{\Delta^n (f(z))}+O\big(\eta^{n+N+1}r^{(n+N+1)(\sigma-1)+\varepsilon}\big)\notag\\
	&= n!\sum_{k=1}^N\frac{\eta^k}{k!}\frac{\Delta^n (f^{(k)})(z)}{f(z)}+O\big(\eta^{n+N+1}r^{(n+N+1)(\sigma-1)+\varepsilon}\big)\notag\\
	&= n!\sum_{k=1}^N\frac{\eta^k}{k!}\frac{f^{(k)}(z)}{f(z)}\frac{\Delta^n (f^{(k)})(z)}{f^{(k)}(z)}+O\big(\eta^{n+N+1}r^{(n+N+1)(\sigma-1)+\varepsilon}\big)\notag\\
		&= n!\sum_{k=1}^N\frac{\eta^k}{k!}\frac{f^{(k)}(z)}{f(z)}\Big(\sum_{s=n}^N\frac{\eta^s\mathfrak{S}_s^{(n)}}{s!}\frac{f^{(k+s)}(z)}{f^{(k)}(z)}+O\big(\eta^{n+N+1}r^{(n+N+1)(\sigma-1)+\varepsilon}\big)\Big)\notag\\
	&\hskip5cm +O\big(\eta^{n+N+1}r^{(n+N+1)(\sigma-1)+\varepsilon}\big)\notag\\
	&=n!\sum_{k=1}^N\frac{\eta^k}{k!}\Big(\sum_{s=n}^N\frac{\eta^s\mathfrak{S}_s^{(n)}}{s!}\frac{f^{(k+s)}(z)}{f(z)}+O\big(\eta^{n+N+1}r^{(n+N+1)(\sigma-1)+\varepsilon}\big)\Big)\notag\\
	&\hskip5cm +O\big(\eta^{n+N+1}r^{(n+N+1)(\sigma-1)+\varepsilon}\big)\notag\\
	&=n!\sum_{k+s=t=n+1}^{2N}\Big(\sum_{k=1}^t\frac{\mathfrak{S}_{t-k}^{(n)}}{k!(t-k)!}\Big)\,\eta^t\,\frac{f^{(t)}(z)}{f(z)}+O\big(\eta^{n+N+1}r^{(n+N+1)(\sigma-1)+\varepsilon}\big)\notag\\
	&= n!\sum_{k+s=t=n+1}^{2N}\Big(\frac{1}{t!}\sum_{k=1}^t\frac{t!}{k!(t-k)!}\,\mathfrak{S}_{k}^{(1)}\mathfrak{S}_{t-k}^{(n)}\Big)\,\eta^t\,\frac{f^{(t)}(z)}{f(z)}+O\big(\eta^{n+N+1}r^{(n+N+1)(\sigma-1)+\varepsilon}\big)\notag\\
	&= n!\sum_{k+s=t=n+1}^{2N}\bigg(\frac{1}{t!}\sum_{k=1}^t\binom{t}{k}\mathfrak{S}_{t-k}^{(1)}\mathfrak{S}_{k}^{(n)}\bigg)\,\eta^t\,\frac{f^{(t)}(z)}{f(z)}+O\big(\eta^{n+N+1}r^{(n+N+1)(\sigma-1)+\varepsilon}\big)\label{E:apply-stirling}\\
	&= n!\sum_{t=n+1}^{2N}\frac{1}{t!}\binom{n+1}{1}\mathfrak{S}_{t}^{(n+1)}\eta^t\frac{f^{(t)}(z)}{f(z)}+O\big(\eta^{n+N+1}r^{(n+N+1)(\sigma-1)+\varepsilon}\big)\notag\\
	&= (n+1)!\sum_{t=n+1}^{2N}\frac{\eta^t\,\mathfrak{S}_{t}^{(n+1)}}{t!}\frac{f^{(t)}(z)}{f(z)}+O\big(\eta^{n+N+1}r^{(n+N+1)(\sigma-1)+\varepsilon}\big)\notag\\
	&= (n+1)!\sum_{t=n+1}^{N}\frac{\eta^t\,\mathfrak{S}_{t}^{(n+1)}}{t!}\frac{f^{(t)}(z)}{f(z)}+O\big(\eta^{n+N+1}r^{(n+N+1)(\sigma-1)+\varepsilon}\big),\notag
\end{align}
where (\ref{E:apply-stirling}) follows from applying $r=1$ from the recurrence formula (\ref{E:recurrence}). 
Each manipulation of above argument is valid outside an exceptional set of $|z|$ of finite logarithmic measure. Since there is at most a finite union of such exceptional sets involved,  so that the above argument is valid outside an exceptional set of finite logarithmic measure. This completes the proof.
\end{proof}

\section{Sharp Difference Wiman-Valiron Estimates}\label{S:Wiman-Valiron}

Let $f(z)=\sum_{n=0}^{\infty}a_nz^n$ be an entire function defined in the complex plane. We recall that  $M(r, f)=\max_{|z|=r}|f(z)|$ denotes the maximum modulus of $f$ on $|z|=r>0$, and we use $\mu(r, f)=\max_{n\geq 0}|a_n|r^n$ to denote the \textit{maximal term} of $f$. The \textit{central index} $\nu(r, f)$ is the greatest exponent $m$ such that
\begin{equation*}
    |a_m|r^m=\mu(r, f).
\end{equation*}
We note that $\nu(r, f)$ is a real, non-decreasing function of
$r$.
\medskip

It is well-known that for finite order $\sigma$ function $f$ that its central
index satisfies $\limsup_{r\rightarrow\infty}{\log\nu(r, f)}/{\log r}
     =\sigma$ (see \cite{Jank_Volkmann1985}, \cite{Laine} and \cite{Wittich}, and also \cite{Fenton2006}). We next quote the classical result of Wiman-Valiron (see also \cite{Hay73}) in the form

\begin{lemma}
\label{L:WM-lemma-2}\cite[pp.
28--30]{HS1988} Let $f$ be a transcendental entire function. Let
$0<~\varepsilon<\frac18$ and $z\ (|z|=r)$ be such that 
\begin{equation}
\label{E:Wiman-Valiron-1}
     |f(z)|>M(r, f)(\nu(r, f))^{-\frac18+\varepsilon}
\end{equation}
holds. Then there exists a set $E\subset (1,\infty)$ of finite
logarithmic measure, such that
\begin{equation}
\label{E:Wiman-Valiron-2}
    \frac{f^{(k)}(z)}{f(z)}=\Big(\frac{\nu(r,
    f)}{z}\Big)^k(1+R_k(z)),
\end{equation}
\begin{equation}
\label{E:Wiman-Valiron-3}
    R_k(z)=O((\nu(r,f))^{-\frac18+\varepsilon})
\end{equation}
holds for all  $k\in \mathbb{N}$ and all $r\notin E\cup [0,1]$.
\end{lemma}

Suppose in addition that we assume $\sigma(f)=\sigma<1$, 
$0<\varepsilon<\frac18$ hold, and that $|z|=r$ is chosen that satisfies (\ref{E:Wiman-Valiron-1}). Then, we deduce, as we have done in \cite{CF-2}, that from (\ref{E:delta-deriv-1}),
(\ref{E:Wiman-Valiron-2}) and (\ref{E:Wiman-Valiron-3}) that for each positive integer $k$, there exists a set
$E\subset (1,\infty)$ of finite logarithmic measure, such
that for all $r\notin E\cup [0,1]$,
\begin{equation*}
\label{E:WV-zero-order}
    \frac{\Delta^kf(z)}{f(z)}=\Big(\frac{\nu(r, f)}{z}\Big)^k
    \big(1+O((\nu(r, f))^{-\frac18+\varepsilon})\big),\quad \textrm{if}\quad \sigma=0,
\end{equation*}
\begin{equation}
\label{E:WV-pos-order}
    \frac{\Delta^kf(z)}{f(z)}=\Big(\frac{\nu(r, f)}{z}\Big)^k
    +O(r^{k\sigma-k-\gamma+\varepsilon}), \quad \textrm{if}\quad 0<\sigma<1,
\end{equation}
where $\gamma=\min\{\frac18\sigma, 1-\sigma\}$.
\bigskip

The remainder in (\ref{E:WV-pos-order}) is not sharp for entire functions of order less than one but not completely regular growth. Equation \eqref{E:higher-difference} from 
Theorem \ref{T:higer-difference}  allows us to remove this restriction and to establish a sharp error bound on (\ref{E:WV-pos-order}).

\begin{theorem}
\label{T:Sharp-WV} Let $f$ be a transcendental entire
function of order $\sigma(f)=\sigma<1$, 
$0<\varepsilon<\min\{\frac18,\, 1-\sigma\}$ and $z$  satisfies (\ref{E:Wiman-Valiron-1}). Then for each positive integer $k$, there exists a set
$E\subset (1,\infty)$ that has finite logarithmic measure, such
that for all $r\notin E\cup [0,\,1]$,
\begin{equation}
\label{E:Sharp-WV}
    \frac{\Delta^kf(z)}{f(z)}=\Big(\frac{\nu(r, f)}{z}\Big)^k \big(1+\mathcal{R}_k(z)\big)
\end{equation}
where $\mathcal{R}_k(z)=O\big(\nu(r,\, f)^{-\kappa+\varepsilon}\big)$ and $\kappa=\min\{\frac18, 1-\sigma\}$.
\end{theorem} 

\medskip 

\begin{proof} Let $0<\varepsilon<\min\{\frac18,\, 1-\sigma\}$ be given. We choose $N$ so large such that
\begin{equation}\label{E:length}
	(n+N+1)(\sigma-1)< -n-1/8.
\end{equation}
We also need
	\begin{equation}\label{E:Wiman}
		\nu(r,f)=O(r^{\sigma+\varepsilon})=O(r),\quad r^{-1}
		=O\big( \nu(r,f)^{\frac{-1}{\sigma+\varepsilon}}\big)=O\big(\nu(r,f)^{-1}\big).
	\end{equation}
Substituting (\ref{E:Wiman-Valiron-2}) into (\ref{E:higher-difference}), applying \eqref{E:length} and \eqref{E:Wiman} yield
\begin{align*}
	\frac{\Delta^n f}{f} &=n!\, \Bigg(\sum_{k=n}^N\frac{\eta^k\,\mathfrak{S}_k^{(n)}}{k!}\frac{f^{(k)}(z)}{f(z)}\Bigg)+O\big(\eta^{n+N+1}r^{(n+N+1)(\sigma-1)+\varepsilon}\big)\notag\\
	&=n!\sum_{k=n}^N\frac{\eta^k\mathfrak{S}_k^{(n)}}{k!}\Big(\frac{\nu(r,
    f)}{z}\Big)^k\big(1+R_k(z)\big)+O\big(\eta^{n+N+1}r^{(n+N+1)(\sigma-1)+\varepsilon}\big)\notag\\
    &= \Big(\frac{\nu(r,f)}{z}\Big)^n\Big[1+\sum_{k=n+1}^Nn!\,\frac{\eta^k\mathfrak{S}_k^{(n)}}{k!}\Big(\frac{\nu(r, f)}{z}\Big)^{k-n}\Big]\big(1+O(\nu (r)^{-\frac{1}{8}+\varepsilon} )\big)\\
	&\qquad +O\big(r^{(n+N+1)(\sigma-1)+\varepsilon}\big)\notag\\
    &= \Big(\frac{\nu(r,f)}{z}\Big)^n\Big[1
+O\big(\nu (r\, f)^{-\frac{1}{8}+\varepsilon}\big)+ O\Big(\frac{\nu(r,\,f)}{r}\Big)\Big]+O(r^{-n-1/8+\varepsilon})\\
	&=\Big(\frac{\nu(r,f)}{z}\Big)^n\Big[1
+O\big(\nu (r\, f)^{-\frac{1}{8}+\varepsilon}\big)+ O\big(\nu(r,\,f)^{1-\frac{1}{\sigma+\varepsilon}}\big)\Big]\\
&=\Big(\frac{\nu(r,f)}{z}\Big)^n\Big[1
+O\big(\nu (r\, f)^{-\frac{1}{8}+\varepsilon}\big)+ O\big(\nu(r,\,f)^{\sigma-1+\varepsilon}\big)\Big]
\end{align*}
as required. This completes the proof.
\end{proof}
\smallskip

We now show that the estimate \eqref{E:Sharp-WV} is sharp in the sense that it is no longer valid for order one entire functions in general. Let $\Phi(z)=1/\Gamma(z)$ where the $\Gamma(z)$ denotes Euler-gamma function. Then it follows from the functional equation $\Gamma(z+1)=z\Gamma(z)$ that
	\[
		\frac{\Delta \Phi(z)}{\Phi(z)}=\frac1z-1\sim -1
	\]
\medskip

\noindent as $z\to\infty$. However, the function $\Phi(z)$ has order equal to one.
\bigskip

\section{Applications to Difference Equations}\label{S:application}
We considered linear difference equations
\begin{equation}
\label{E:difference-eq-1}
    a_n(z)\Delta^n f(z)+\cdots+a_{1}(z)\Delta f(z)+a_0(z)f(z)=0
\end{equation}
in \cite{CF-2}, where  $a_0(z),\,\cdots, a_n(z)$ are polynomials. We have shown that any entire solution $f$ to (\ref{E:difference-eq-1}) with order of growth less than one has a positive rational order of growth and that the rational order $\chi$ can be calculated explicitly from the gradients of the
corresponding Newton-Puiseux diagram of the equation.  
This falls short of showing these entire solutions have completely regular growth as is well-known for entire solutions to linear differential equations (\textit{i.e.,} replacing the $\Delta^k f$ in (\ref{E:difference-eq-1}) by $f^{(k)}$ for $k=1,\cdots, n$). We shall strengthen our earlier result that the growth order $\chi$ of entire solutions to (\ref{E:difference-eq-1}) with $\chi<1$ is indeed completely regular. This also improves an earlier result of Ishizaki and Yanagihara \cite{IY2004} where they have proved the following result for solutions under the assumption of growth order $<1/2$ by developing a Wiman-Valiron theory based on \textit{binomial series}. We mention in passing that Ramis considered the corresponding problems for $q-$difference equations \cite{Ram}.
\bigskip

\begin{theorem}
\label{T:application-1} Let $a_0(z),\,\cdots, a_n(z)$ be polynomial
coefficients of the difference equation (\ref{E:difference-eq-1}), and let $f$ be an entire solution with order
$\sigma (f)=\chi< 1$. Then $\chi$
is a rational number which can be determined from a gradient of the
corresponding Newton-Puisseux diagram equation
(\ref{E:difference-eq-1}). In particular,
\[
	\log M(r,\, f)=Lr^\chi\big(1+o(1)\big)
\]
where $L>0$, $\chi>0$ and $M(r,\, f)=\max_{|z|=r}|f(z)|$. That is, the solution has completely regular growth.
\end{theorem}

The original proof of the corresponding result for linear differential equations is based on an application of (\ref{E:Wiman-Valiron-2}). Thus the form of our estimate (\ref{E:Sharp-WV}) allows us to modify Valiron's argument \cite{Val}
simply by replacing the differential operators with difference operators and the rest of the proof follows exactly the same pattern. So we omit the details. Wittich \cite[pp. 65--68]{Wittich} and Gundersen, Steinbart and Wang \cite{GSW98} discussed if the possible rational order for the entire solutions of linear differential equations with polynomial coefficients obtained by the Newton-Puiseux diagram (\cite[\S 22]{Jank_Volkmann1985})
method could be attained. It would be interesting to ask this question for difference equation (\ref{E:difference-eq-1}).
\bigskip

\section*{Acknowledgment}
The authors would like to thank their colleague Dr. T. K. Lam for pointing out an inaccuracy on the exceptional set, now labelled as $E^{(\eta)}$, in an earlier version of this paper. The first author would also like to acknowledge hospitality that he received during his visits to Academy of Mathematics and Systems Science, Chinese Academy of Sciences, Beijing in the preparation of this paper.  Finally, both authors thank the referees for their constructive and detailed comments that helped to improve the readability of this paper.
\medskip


\begin{thebibliography}{}
%
%
\bibitem{AHH} {M. J. Ablowitz, R. Halburd and B. Herbst}, \textrm{On the extension of the Painlev\'e property to difference  equations}, \textit{Nonlinearity} \textbf{13}, (2000), no. 3, 889--905.

\bibitem{AS64} M. Abramowitz and I. A. Stegun (ed.), \textit{Handbook of Mathematical Functions}, National Bureau of Standards, Applied Mathematics Series 55, Washington, 1964.

\bibitem{Apo} {T. M. Apostol}, \textrm{Mathematical Analysis}, Wiley 1967.

\bibitem{BL2007} {W. Bergweiler and J. K. Langley}, \textrm{``Zeros of differences of meromorphic functions"}, \textit{Math. Proc. Cambridge Philos. Soc.} \textbf{142} (2007) no. 1, 133--147.

\bibitem{Cartan1928} {H. Cartan}, \textrm{``Sur les syst\`ems de fonctions holomorphes \`a vari\'et\'es lin\'eaires lacunaires et leurs applications"}, \textit{Ann. Sci. Ecole Norm. Sup.} (3) \textbf{45} (1928) 255--346.

\bibitem{CF-1} {Y. M. Chiang and S. J. Feng}, \textrm{``On the
Nevanlinna characteristic of $f(z+\eta)$ and difference equations
in the complex plane"}, \textit{The Ramanujan J.},  \textbf{16} (2008), no. 1, 105--129.

\bibitem{CF-2} {Y. M. Chiang and S. J. Feng}, \textrm{``On the growth
of logarithmic differences, difference quotients and logarithmic
derivatives of meromorphic functions"},  \textit{Trans. Amer. Math. Soc.}  \textbf{361}  (2009), no. 7, 3767--3791.



\bibitem{Fenton2006} P. C. Fenton, {``A glance at Wiman-Valiron theory"}, Complex analysis and dynamical systems II, 131--139,  Contemp. Math., 382, Amer. Math. Soc., Providence, RI, 2005. 


\bibitem{Gund88} {G. G. Gundersen}, \textrm{``Estimates for the
logarithmic derivative of meromorphic functions, plus similar
estimates"}, \textit{J. London Math. Soc.} (2), \textbf{37}, (1988), no. 1,  88--104.


\bibitem{GSW98} {G. G. Gundersen, E. M. Steinbart and S. P. Wang},
\textrm{``The possible orders of solutions of linear differential
equations with polynomial coefficients"}, \textit{Trans. Amer. Math. Soc.},
\textbf{350} (1998), no. 3, 1225--1247.

\bibitem{HK-1} {R. G. Halburd and R. J. Korhonen},
\textrm{``Difference analogue of the lemma on the logarithmic
derivative with applications to difference equations"}, \textit{J. Math.
Anal. Appl.} \textbf{314}, (2006), no. 2, 477--487.

\bibitem{HK-2} {R. G. Halburd and R. J. Korhonen},
\textrm{``Finite-order meromorphic solutions and the discrete
Painlev\'e equations"}, \textit{Proc. London Math. Soc.}, (3) \textbf{94} (2007). no. 2, 443--474.


\bibitem{HK-5} {R. G. Halburd and R. J. Korhonen}, \textrm{``Meromorphic solutions of difference equations, integrability and the discrete Painlev\'e equations"}, \textit{J. Phys. A} \textbf{40} (2007), no. 6, R1--R38.
 
\bibitem{Hay64} {W. K. Hayman}, \textit{Meromorphic Functions}, Oxford Univ.
  Press, Oxford, 1964.

\bibitem{Hay73} {W. K. Hayman}, \textrm{``The local growth of power series
: a survey of the Wiman-Valiron method"}, \textit{Canad. Math. Bull.}
\textbf{17} (1974), no. 3, 317--358.


\bibitem{HS1988} {Y. He and X. Xiao}, \textit{Algebriod Functions
and Ordinary Differential Equations}, Beijing Sci. Press, 1988
(Chinese).

\bibitem{HKLRT} {J. Heittokangas, R. Korhonen, I. Laine, J.
Rieppo and K. Tohge}, \textrm{``Complex difference equations of
Malmquist type"}, \textit{Comp. Methods Funct. Theory}, \textbf{1} (2001), no. 1, 27--39.


\bibitem{IY2004} {K. Ishizaki and N. Yanagihara},
\textrm{``Wiman-Valiron method for difference equations"}, \textit{Nagoya Math.
J.}, \textbf{175} (2004), 75--102.

\bibitem{Jank_Volkmann1985} {G. Jank and L.Volkmann}, \textit{Einf\"uhrung in die Theorie der ganzen und meromorphen Funktionen mit Anwendungen auf Differentialgleichungen}, UTB f\"ur Wissenschaft: Gross Reihe. 
 Birkh\"auser Verlag, Basel, 1985.

\bibitem{Laine} {I. Laine}, \textit{Nevanlinna Theory and Complex
Differential Equations}, de Gruyter Studies in Mathematics, 15. Walter de Gruyter \& Co., Berlin 1993.



\bibitem{Ram} {J.-P. Ramis}, \textrm{``About the growth of entire
function solutions of linear algebraic $q$--difference equations"},
 \textit{Ann. Fac. Sci. Toulouse Math.}, (6)  \textbf{1}  (1992),  no. 1, 53--94.


\bibitem{Val} {G. Valiron}, \textit{Lectures on the Theory of Integral
Functions}, Chelsea, Pulb. Co., 1949 (translated by E. F. Collingwood).

\bibitem{Wittich} Hans Wittich, \textit{Neuere Untersuchungen \"uber eindeutige analytische Funkitonen} (2nd ed.)
 Ergebnisse d. Math, N. S., No. 8. Springer-Verlag, Berlin, 1955; 2nd ed., 1968

\end{thebibliography}
\end{document}